\documentclass[12pt,thmsa]{article}
\usepackage{amsfonts}

\newtheorem{theorem}{Theorem}
\newtheorem{corollary}{Corollary}
\newtheorem{proposition}{Proposition}
\newtheorem{lemma}{Lemma}
\newtheorem{rem}{Remark}
\newcommand{\p}{\Bbb{P}}

\newcommand{\pf}{\Bbb{P}^\uparrow}
\newcommand{\pfx}{\Bbb{P}_{x}^\uparrow}

\newcommand{\ind}{\mbox{\rm 1\hspace{-0.04in}I}}

\newcommand{\eqdef}{\stackrel{\mbox{\tiny$($def$)$}}{=}}

\def\QED{\hfill\vrule height 1.5ex width 1.4ex depth -.1ex \vskip20pt}

\begin{document}

\vspace*{-0.5in}
 \hspace*{-0.5in} {\footnotesize This version: 31/05/2005}

\vspace*{0.9in}
\begin{center}
{\LARGE  On L\'{e}vy processes conditioned to stay positive.\vspace*{0.4in}}\\
{\large L. Chaumont$^*$ and R.A. Doney$^\sharp$\vspace*{0.2in}}
\end{center}
\noindent {\footnotesize $^*$Laboratoire de Probabilit\'es et
Mod\`eles Al\'eatoires, Universit\'e Pierre et Marie Curie,
4, Place Jussieu - 75252 {\sc Paris Cedex 05.} E-mail: chaumont@ccr.jussieu.fr\\
$^\sharp$Department of Mathematics, University
of Manchester, {\sc Manchester, M 13 9 PL.}\\ E-mail: rad@ma.man.ac.uk}\\

\noindent {\it Abstract}: {\footnotesize We construct the law of
L\'evy processes conditioned to stay positive under general
hypotheses. We obtain a Williams type path decomposition at the
minimum of these processes. This result is then applied to prove
the weak convergence of the law of L\'evy processes conditioned to
stay positive as their initial state tends to 0. We describe an
absolute continuity relationship between the limit law and the
measure of the excursions away from 0 of the underlying L\'evy
process reflected at its minimum. Then, when the L\'evy process
creeps upwards, we study the lower tail at
0 of the law of the height this excursion.}\\

\noindent {\it Key words}: {\footnotesize  L\'evy process
conditioned to stay positive, path decomposition, weak
convergence, excursion measure, creeping.}

\noindent
{\it A.M.S. Classification}: {\footnotesize 60G51, 60G17.}\\

\section{Introduction}
\setcounter{equation}{0}

In \cite{c1} it was shown how, given the measure $\Bbb{P}$ of a
L\'{e}vy process satisfying some weak assumptions, one can
construct for each $x>0$ a measure $\Bbb{P}_{x}^{\uparrow }$
corresponding, in the sense of Doob's theory of h-transforms, to
conditioning the process starting at $x$ to stay positive. Using a
different construction Bertoin \cite{b3} had shown the existence
of a measure $\Bbb{P}_{0}^{\uparrow }$ under which the process
starts at $0$ and stays positive. A natural question is whether
$\Bbb{P} _{x}^{\uparrow }$ converges to $\Bbb{P}^{\uparrow
}:=\Bbb{P}_{0}^{\uparrow }$ as $x\downarrow 0.$ Recently, it has
been proved by Tanaka \cite{ta} that this convergence in law holds
in the sense of finite dimensional distributions, under very
general hypotheses, but here we are interested in convergence in
law on Skohorod's space of c\`{a}dl\`{a}g trajectories. Since it
was also shown in \cite{c1}, extending a famous result for the
3-dimensional Bessel process due to Williams, that under $\Bbb{P}
_{x}^{\uparrow }$ the post-minimum process is independent of the
pre-minimum process and has law $\Bbb{P}^{\uparrow }$, this
essentially amounts to showing that the pre-minimum process
vanishes as $x\downarrow 0.$ Such a result has been verified in
the case of spectrally negative processes in \cite{b2}, and for
stable processes and for processes which creep downwards in
\cite{c1}.\newline

In the third section of this paper, we give a simple proof of this
result for a general L\'{e}vy process. This proof does not use the
description of the law of the pre-minimum process which is given
in \cite{c1}, but is based on knowledge of the distribution of the
all-time minimum under $\Bbb{P}_{x}^{\uparrow },$ which was also
established in \cite{c1}. As a consequence, we are able to extend
the description of the excursion measure of the process reflected
at its minimum (see Proposition 15, P. 202 of \cite {b1} for the
spectrally negative case) to the general case. This description
has recently been used in \cite{d2} to perform some semi-explicit
calculations for the reflected process when $X$ has jumps of one
sign only, and our result suggests the possibility that such
calculations could be performed in other cases. However a key fact
in the calculations in \cite{d2} is that, in the spectrally
negative case, the harmonic function $h$ used in conditioning to
stay positive and the excursion measure $\underline{n}$ of the
reflected process satisfy
\[
\underline{n}(H>x)h(x)=1,
\]
where $H$ denotes the height of a generic excursion. This leads us
to wonder whether this relation could hold \textbf{asymptotically
as }$x\downarrow 0$ in other cases. Intuitively it seems obvious
that this should be the case when the reflected excursion
\textbf{creeps upwards}. This is what we establish in Theorem
\ref{theo123}.

\section{Notation and Preliminaries}
\setcounter{equation}{0}

Let $\mathcal{D}$ denote the space of c\`{a}dl\`{a}g paths $\omega
:[0,\infty )\rightarrow \Bbb{R\cup \{\delta \}}$ with lifetime
$\zeta (\omega )=\inf \{s:\omega _{s}=\delta \},$ where $\delta $
is a cemetery point. $\mathcal{D}$ will be equipped with the
Skorokhod topology, with its Borel $\sigma $-algebra
$\mathcal{F}$, and the usual filtration $(\mathcal{F} _{s},s\geq
0).$ Let $X=(X_{t},t\geq 0)$ be the coordinate process defined on
the space $\mathcal{D}$. We write $\overline{X}$ and
$\underline{X}$ for the supremum and infimum processes, defined
for all $t<\zeta $ by
\begin{eqnarray*}
\overline{X}_{t} &=&\sup \{X_{s}:0\leq s\leq t\}, \\
\underline{X}_{t} &=&\inf \{X_{s}:0\leq s\leq t\}.
\end{eqnarray*}
We write $\tau _{A}$ for the entrance time into a Borel set $A,$
and $m$ for the time at which the absolute infimum is attained:
\begin{eqnarray*}
\tau _{A} &=&\inf \{s>0:X_{s}\in A\}, \\
m &=&\sup \{s<\zeta :X_s\wedge X_{s-}=\underline{X}_{s}\},
\end{eqnarray*}
with the usual convention that $\inf \{\emptyset \}=+\infty \mbox{
and }\sup \{\emptyset \}=0$.

For each $x\in \Bbb{R}$ we denote by $\Bbb{P}_{x}$ the law of a
L\'{e}vy process starting from $x,$ and write
$\Bbb{P}_{0}=\Bbb{P}$. We assume throughout the sequel that
$(X,\Bbb{P})$ is not a compound Poisson process. It is well known
that the reflected process $X-\underline{X}$ is Markov. Note that
the state 0 is regular for $(-\infty,0)$ under $\Bbb{P}$, if and
only if it is regular for $\{0\}$ for the reflected process. In
this case, we will simply say that 0 is regular downwards and if 0
is regular for $ (0,\infty)$ under $\Bbb{P}$, we will say that 0
is regular upwards.

Let $\underline{L}$ be the local time of the reflected process
$X-\underline{ X}$ at 0 and let $\underline{n}$ be the measure of
its excursions away from 0. If 0 is regular downwards then, up to
a multiplicative constant, $ \underline{L}$ is the unique additive
functional of the reflected process whose set of increasing points
is $\{t:(X-\underline{X})_t=0\}$ and $ \underline{n}$ is the
corresponding It\^o measure of excursions; we refer to \cite{b1},
Chap. IV, sections 2-4 for a proper definition of $\underline{L}$
and $\underline{n}$. If 0 is not regular downwards then the set
$\{t:(X- \underline{X})_t=0\}$ is discrete and we define the local
time $\underline{L} $ as the counting process of this set, i.e.
$\underline{L}$ is a jump process whose jumps have size 1 and
occur at each zero of $X-\underline{X}$. Then, the measure
$\underline{n}$ is the probability law of the process $X$ under
the law $\Bbb{P}$, killed at its first passage time in the
negative halfline, i.e. $\tau_{(-\infty,0)}$, (see the definition
of $\Bbb{Q}_0$ below).

Let us first consider the function $h$ defined for all $x\ge 0$
by:
\begin{equation}  \label{h}
h(x):=\Bbb{E}\left(\int_{[0,\infty )}1_{\{\underline{X}_{t}\geq
-x\}}\,d \underline{L}_{t}\right) \,. \label{eq:id1}
\end{equation}
It follows from (\ref{h}) (or (\ref{lim}) below) and general
properties of L\'{e}vy processes that $h$ is {\it finite,
continuous, increasing} and that $h-h(0)$ is {\it subadditive} on
$[0,\infty )$. Moreover, $h(0)=0$ if 0 is regular downwards and
$h(0)=1$ if not (in the latter case, the counting measure
$d\underline{L}_{t} $ gives mass 1 to the point $t=0$).

Let ${\bf e}$ be an exponential time with parameter 1, which is
independent of $(X,\Bbb{P})$. The following identity follows from
the exit formula of excursion theory when 0 is regular downwards
and is obtained by direct calculations in the other case. For all
$\varepsilon>0$,
\begin{equation}\label{exit}
\Bbb{P}_{x}(\tau _{(-\infty ,0)}>\mbox{{\bf e}}/\varepsilon
)=\Bbb{P}( \underline{X}_{\mbox{{\bf e}}/\varepsilon }\geq
-x)=\Bbb{E}\left( \int_{[0,\infty)}e^{-\epsilon
t}\ind_{\{\underline{X}_{t}\geq -x\}}\,d \underline{L}_{t}\right)
\,\underline{n}(\mbox{{\bf e}}/\varepsilon <\zeta )\,,
\end{equation}
so that, by monotone convergence, for all $x\ge 0$:
\begin{equation}  \label{lim}
h(x)=\lim_{\varepsilon \rightarrow 0}\frac{\Bbb{P}_{x}(\tau
_{(-\infty ,0)}> \mbox{{\bf
e}}/\varepsilon)}{\underline{n}(\mbox{{\bf e}}/\varepsilon
<\zeta)}\,.
\end{equation}
In the next lemma, we show that $h$ is excessive or invariant for
the process $(X,\Bbb{P}_{x})$, $x>0$ killed at time $\tau
_{(-\infty ,0)}$. This result has been proved by Silverstein
\cite{s}, see also Tanaka \cite{ta}. Here, we give a different
proof which uses the representation of $h$ stated in (\ref{lim}).
For $x>0$, we denote by $\Bbb{Q}_{x}$ the law of the killed
process, i.e. for $\Lambda \in \mathcal{F}_{t}$:
\[
\Bbb{Q}_{x}(\Lambda ,t<\zeta )=\Bbb{P}_{x}(\Lambda ,\,t< \tau
_{(-\infty ,0)})\,,
\]
and by $(q_{t})$ its semigroup. Recall that $\Bbb{Q}_{0}$ and
$q_t(0,dy)$ are well defined when 0 is not regular downwards, and
in this case we have $\Bbb{Q}_{0}=\underline{n}$.

\begin{lemma}
\label{lem12} If $(X,\Bbb{P})$ drifts towards $-\infty $ then $h$
is excessive for $(q_{t})$, i.e. for all $x\ge0$ and $t\ge 0$,
$\Bbb{E }_{x}^{\Bbb{Q}}(h(X_{t})\mbox{\rm
1\hspace{-0.033in}I}_{\{t<\zeta \}})\le h(x)$. If $(X,\Bbb{P})$
does not drifts to $-\infty $, then $h$ is invariant for
$(q_{t})$, i.e. for all $x\ge0$ and $t\ge 0$,
$\Bbb{E}_{x}^{\Bbb{Q} }(h(X_{t})\mbox{\rm
1\hspace{-0.033in}I}_{\{t<\zeta \}})=h(x)$.
\end{lemma}

\noindent {\it Proof}: From (\ref{lim}), monotone convergence and
the Markov property, we have
\begin{eqnarray}
\Bbb{E}^{\Bbb{Q}}_x(h(X_t)\mbox{\rm
1\hspace{-0.033in}I}_{\{t<\zeta\}})&=&
\lim_{\varepsilon\rightarrow0} \Bbb{E}_x\left(\frac{\Bbb{P}
_{X_t}(\tau_{(-\infty,0)}>\mbox{{\bf e}}/\varepsilon)\mbox{\rm
1\hspace{-0.033in}I}_{\{t\le\tau_{(-\infty,0)}\}}}
{\underline{n}(\mbox{{\bf
e}}/\varepsilon<\zeta)}\right)  \nonumber \\
&=&\lim_{\varepsilon\rightarrow0} \Bbb{E}_x\left(\frac{\mbox{\rm
1\hspace{-0.033in}I}_{\{\tau_{(-\infty,0)}>t+\mbox{{\bf
e}}/\varepsilon\}}}{ \underline{n}(\mbox{{\bf
e}}/\varepsilon<\zeta)}\right) =\lim_{\varepsilon\rightarrow0}
e^{\varepsilon t}\left(\frac{\Bbb{P}
_{x}(\tau_{(-\infty,0)}>\mbox{{\bf
e}}/\varepsilon)}{\underline{n}(
\mbox{{\bf e}}/\varepsilon<\zeta)}\right.  \nonumber \\
&&\left.-\int_0^t\varepsilon e^{-\varepsilon u}\frac{\Bbb{P}
_{x}(\tau_{(-\infty,0)}>u)}{\underline{n}(\mbox{{\bf
e}}/\varepsilon<\zeta)}
\,du\right)  \nonumber \\
&=&h(x)-\frac{1}{\underline{n}(\zeta)}\int_0^t\Bbb{P}_{x}(\tau_{(-
\infty,0)}>u)\,du\,,\label{ex}
\end{eqnarray}
where
$\underline{n}(\zeta)=\int_0^{\infty}\underline{n}(\zeta>t)\,dt$.
It is known that for $x>0$, $E_x(\tau_{(-\infty,0)})<\infty$ if
and only if $ X$ drifts towards $-\infty$, see \cite{b1}, Prop.
VI.17. Hence, since moreover for $x>0$, $0<h(x)<+\infty$, then
(\ref{lim}) shows that $ \underline{n}(\zeta)<+\infty$ if and only
if $X$ drifts towards $-\infty$. Consequently, from (\ref{ex}), if
$X$ drifts towards $-\infty$, then $\Bbb{E}
_x^{\Bbb{Q}}(h(X_t)\mbox{\rm
1\hspace{-0.033in}I}_{\{t<\zeta\}})\le h(x)$, for all $t\ge0$ and
$x\ge0$, whereas if $(X,\Bbb{P})$ does not drift to $-\infty$,
then $\underline{n} (\zeta)=+\infty$ and (\ref{ex}) shows that
$\Bbb{E}_x^{\Bbb{Q}}(h(X_t) \mbox{\rm
1\hspace{-0.033in}I}_{\{t<\zeta\}})=h(x)$, for all $t\ge0$ and
$x\ge0$.\QED

\section{The process conditioned to stay positive}

\setcounter{equation}{0}

We now define the L\'evy process $(X,\Bbb{P}_x)$ conditioned to
stay positive. This notion  has now a long history, see \cite{b3},
\cite{c1} , \cite{du}, \cite{ta} and the references contained in
these papers. We begin this section with some properties of the
process conditioned to stay positive (stated in Proposition
\ref{prop1} and Theorem \ref{theo1}) which are usually established
under the additional assumptions that the semigroup of
$(X,\Bbb{P})$ is absolutely continuous and/or 0 is regular
downwards and/or upwards...

Write $(p_{t},t\geq 0)$ for the semigroup of $(X,\Bbb{P})$ and
recall that $(q_{t},t\geq 0)$ is the semigroup (in $(0,\infty )$
or in $[0,\infty)$) of the process $(X,\Bbb{Q}_{x})$. Then we
introduce the new semigroup
\begin{equation}
p_{t}^{\uparrow }(x,dy):=\frac{h(y)}{h(x)}q_{t}(x,dy),\mbox{
}x>0,y>0,\mbox{ }t\geq 0\,.  \label{eq1}
\end{equation}
From Lemma \ref{lem12}, $(p_{t}^{\uparrow })$ is sub-Markov when
$(X,\Bbb{P}) $ drifts towards $-\infty $ and it is Markov in the
other cases. For $x>0$ we denote by $\Bbb{P}_{x}^{\uparrow }$ the
law of the strong Markov process started at $x$ and whose
semigroup in $(0,\infty )$ is $(p_{t}^{\uparrow })$ . When
$(p_{t}^{\uparrow })$ is sub-Markov, $(X,\Bbb{P}_{x}^{\uparrow })$
has state space $(0,\infty )\cup \{\delta \}$ and this process has
finite lifetime. In any case, for $\Lambda \in \mathcal{F}_{t}$,
we have
\begin{equation}
\Bbb{P}_{x}^{\uparrow }(\Lambda ,t<\zeta
)=\frac{1}{h(x)}\Bbb{E}_{x}^{\Bbb{Q} }(h(X_{t})\mbox{\rm
1\hspace{-0.033in}I}_{\Lambda }\mbox{\rm
1\hspace{-0.033in}I}_{\{t<\zeta \}})\,.  \label{pfx}
\end{equation}
Note that when 0 is not regular downwards then definitions
(\ref{eq1}) and ( \ref{pfx}) also make sense for $x=0$. We show in
the next proposition that $ \Bbb{P}_{x}^{\uparrow }$ is the limit
as $\varepsilon \downarrow 0$ of the law of the process under
$\Bbb{P}_{x}$ conditioned to stay positive up to an independent
exponential time with parameter $\varepsilon ,$ so we will refer
to $(X,\Bbb{P}_{x}^{\uparrow })$ as the process ``conditioned to
stay positive''. Note that the following result has been shown in
\cite{c1} under stronger hypotheses; however, to make this paper
self-contained, we give the proof below.

\begin{proposition}
\label{prop1} Let \textbf{e} be an exponential time with parameter
$1$ which is independent of $(X,\Bbb{P})$.

For any $x>0$, and any $(\mathcal{F}_{t})$ stopping time $T$ and
for all $\Lambda \in \mathcal{F}_{T}$,
\[
\lim_{\varepsilon \rightarrow 0}\Bbb{P}_{x}(\Lambda ,T<\mbox{{\bf
e}} /\varepsilon \,|\,X_{s}>0,\,0\leq s\leq \mbox{{\bf
e}}/\varepsilon )=\Bbb{P} _{x}^{\uparrow }(\Lambda ,T<\zeta )\;.
\]
This result also holds for $x=0$ when $0$ is not regular
downwards.
\end{proposition}

\noindent {\it Proof}: According to the Markov property and the
lack-of-memory property of the exponential law, we have
\begin{eqnarray}
&&\Bbb{P}_x(\Lambda,\,T<\mbox{{\bf e}}/\varepsilon\,|
\,X_s>0,\,0\leq s\leq
\mbox{{\bf e}}/\varepsilon) = \nonumber\\
&&\Bbb{E}_x\left(\ind_\Lambda\ind_{\{T<(\mbox{{\bf
e}}/\varepsilon)\wedge\tau_{(- \infty,0)}\}}
\frac{\Bbb{P}_{X_T}(\tau_{(-\infty,0)}\geq \mbox{{\bf e}}
/\varepsilon)} {\Bbb{P}_x(\tau_{(-\infty,0)}\geq \mbox{{\bf
e}}/\varepsilon)} \right) \;.\label{eqx}
\end{eqnarray}
Let $\varepsilon_0>0$. From (\ref{h}) and (\ref{exit}), for all
$\varepsilon\in(0,\varepsilon_0)$,
\begin{eqnarray}
&&\ind_{\{T<(\mbox{{\bf e}}/\varepsilon)\wedge\tau_{(-
\infty,0)}\}}\frac{\Bbb{P}_{X_T}(\tau_{(-\infty,0)}\geq \mbox{{\bf
e}}/\varepsilon)}{\Bbb{P}_x(\tau_{(-\infty,0)}\geq \mbox{{\bf
e}}/\varepsilon)}\le\nonumber \\ &&\ind_{\{T<\tau_{(-
\infty,0)}\}}\Bbb{E}\left(\int_{[0,\infty)} e^{-\varepsilon_0
t}\ind_{\{\underline{X}_t\ge-x\}}\,d\underline{L}_t\right)^{-1}
h(X_T)\,,\;\;\;\mbox{\rm a.s.}\label{eqy} \end{eqnarray}
 Recall that $h$ is excessive
for the semigroup $(q_t)$, hence the inequality of Lemma
\ref{lem12} also holds at any stopping time, i.e.
$\Bbb{E}^{\Bbb{Q} }_x(h(X_T)\ind_{\{T<\zeta\}})\le h(x)$. Since
$h$ is finite, the expectation of the right hand side of
(\ref{eqy}) is finite so that we may apply Lebesgue's theorem of
the dominated convergence in the right hand side of (\ref{eqx})
when $\varepsilon$ goes to 0. We conclude by using the
representation of $h$ in (\ref{lim}) and the definition of
$\Bbb{P}^\uparrow_x$ in (\ref{pfx}).\QED

\begin{rem} In the discrete time setting, that is for random walks,
another characterization of the harmonic function $h$ has been
given by Bertoin and  Doney {\rm \cite{bd}}, Lemma $1$. Using
similar arguments, one can show that a continuous time equivalent
holds. For L\'evy processes, such that $\limsup_tX_t=+\infty$,
this result is
\[\lim_{n\rightarrow+\infty}\frac{\p_x(\tau_{[n,\infty)}<\tau_{(-\infty,0)})}
{\p_y(\tau_{[n,\infty)}<\tau_{(-\infty,0)})}=\frac{h(x)}{h(y)}\,,\;\;\;x,y>0.\]
Then, as in discrete time, a consequence is the following
equivalent definition of $(X,\pf_x)$$:$
\[\lim_{n\rightarrow+\infty}\p_x(\Lambda\,|\,\tau_{[n,\infty)}<\tau_{(-\infty,0)})=
\pf_x(\Lambda)\,,\;\;t>0,\,\;\;\Lambda\in{\cal F}_t\,.\]  Note
that a similar conditioning has been studied by Hirano {\rm
\cite{hi}} in some special cases.
\end{rem}

In the case where 0 is regular downwards, definition (\ref {eq1})
doesn't make sense for $x=0,$ but in \cite{b3} it was shown that
in any case, the law of the process
\[
((X-\underline{X})_{g_{t}+s},\,s\le t-g_{t}),\;\;\mbox{where}
\;\;\;g_{t}=sup\{s\le t:(X-\underline{X})_{s}=0\}
\]
converges as $t\rightarrow \infty $ to a Markovian law under which
$X$ starts at $0$ and has semigroup $p_{t}^{\uparrow }$.
Similarly, under additional hypotheses, Tanaka \cite{ta}, Th.7
proved that the process
\[
(X-\underline{X})_{b_{\lambda }+s},\,s\le a_{\lambda }-b_{\lambda
}),\;\,\ \mbox{where}\;\;\left\{
\begin{array}{l}
a_{\lambda }=inf\{t:(X-\underline{X})_{t}>\lambda \} \\
b_{\lambda }=sup\{t\le a_{\lambda }:(X-\underline{X})_{t}=0\}
\end{array}
\right.
\]
converges as $\lambda \rightarrow +\infty $ towards the same law.
We will denote this law by $\Bbb{P}^{\uparrow }.$ Thm 3 of
\cite{c1} gives the entrance law of the process
$(X,\Bbb{P}^{\uparrow })$, (see (\ref{eq2}) below). Note that
Doney \cite {do}, extending a discrete time result from Tanaka
\cite{ta1}, obtained a path construction of $(X,\Bbb{P}^{\uparrow
})$. Another path construction of $(X,\Bbb{P}^{\uparrow })$ is
contained in Bertoin \cite{b3}. These two constructions are quite
different from each other but coincide in the Brownian case.
Roughly speaking, we could say that Doney-Tanaka's construction is
based on a rearrangement of the excursions away from 0 of the
L\'{e}vy process reflected at is minimum whereas Bertoin's
construction consists in sticking together the positive excursions
away from 0 of the L\'{e}vy process itself.

The next theorem describes the decomposition of the process
$(X,\Bbb{P}_{x}^\uparrow)$ at the time of its minimum. It is also
well known in the literature under various hypotheses, see
\cite{c1}, \cite{du}. Here, we have tried to state it under the
weakest possible assumptions.

\begin{theorem}
\label{theo1}  Define the pre-minimum and post-minimum processes
respectively as follows:  $(X_{t}\,,\,0\leq t<m)$ and
$(X_{t+m}-U\,,\,0\leq t<\zeta -m)$, where $U:=X_{m}\wedge X_{m-}$.
\begin{itemize}
\item[$1.$]  Under $\Bbb{P}_{x}^{\uparrow }$, $x>0$, the
pre-minimum and post-minimum processes are independent. The
process $(X,\Bbb{P}_{x}^{\uparrow })$ reaches its absolute minimum
$U$ once only and its law is given by:
\begin{equation}
\Bbb{P}_{x}^{\uparrow }(U\ge y)=\frac{h(x-y)}{h(x)}1_{\{y\leq
x\}}\;. \label{min}
\end{equation}

\item[$2.$]  Under $\Bbb{P}_{x}^{\uparrow }$, the law of the
post-minimum process is $\Bbb{P}^{\uparrow }$. In particular, it
is strongly Markov and does not depend on $x$. The semigroup of
$(X,\Bbb{P}^{\uparrow })$ in $(0,\infty )$ is $(p_{t}^{\uparrow
})$. Moreover, $X_{0}=0$, $\Bbb{P}^{\uparrow }$-a.s. if and only
if $0$ is regular upwards.
\end{itemize}
\end{theorem}

\noindent {\it Proof}: Denote by $\Bbb{P} _x^{\mbox{{\bf
e}}/\varepsilon}$ the law of the process $(X,\Bbb{P}_x)$ killed at
time $\mbox{{\bf e}}/\varepsilon$. Since $(X,\Bbb{P})$ is not a
compound Poisson process, it almost surely reaches its minimum at
a unique time on the interval $[0,{\bf e}/\varepsilon]$. Recall
that by a result of Millar \cite{mi}, pre-minimum and post-minimum
processes are independent under $\Bbb{P}_x^{\mbox{{\bf
e}}/\varepsilon}$ for all $\varepsilon>0$. According to
Proposition \ref{prop1}, the same properties hold under
$\Bbb{P}_{x}^\uparrow$. Let $0\leq y\le x$. From Proposition
\ref{prop1} and (\ref {lim}):
\begin{eqnarray*}
\Bbb{P}_{x}^\uparrow(U<
y)&=&\Bbb{P}_{x}^\uparrow(\tau_{[0,y)}<\zeta)=\lim_{
\varepsilon\rightarrow0}\Bbb{P}_x(\tau_{[0,y)}<\mbox{{\bf e}}
/\varepsilon\,|\, \tau_{(-\infty,0)}>\mbox{{\bf e}}/\varepsilon)\\
&=&\lim_{\varepsilon\rightarrow0}\left(1-\frac{\Bbb{P}_x(\tau_{[0,y)}\ge\mbox{{\bf
e}} /\varepsilon,\,\tau_{(-\infty,0)}>\mbox{{\bf e}}
/\varepsilon)}{\Bbb{P}_x(\tau_{(-\infty,0)}>\mbox{{\bf e}}
/\varepsilon)}\right)\\ &=&1-
\lim_{\varepsilon\rightarrow0}\frac{\Bbb{P}_{x-y}(\tau_{(-\infty,0)}\ge
\mbox{{\bf
e}}/\varepsilon)}{\Bbb{P}_x(\tau_{(-\infty,0)}>\mbox{{\bf e}}
/\varepsilon)}=1-\frac{h(x-y)}{h(x)}\,,
\end{eqnarray*}
and the first part of the theorem is proved.

From the independence mentioned above, the law of the post-minimum
process under $\Bbb{P}_x^{\mbox{{\bf e}}/\varepsilon}(\,\cdot
\,|\,U>0)$ is the same as the law of the post-minimum process
under $\Bbb{P }_x^{\mbox{{\bf e}}/\varepsilon}$. Then, from
Proposition \ref {prop1} or from \cite{b3}, Corollary 3.2, the law
of the post-minimum processes under $\pfx$ is the limit of the law
of the post-minimum process under $\Bbb{P} _x^{\mbox{{\bf
e}}/\varepsilon}$, as $\varepsilon\rightarrow0$. But in \cite{b3},
Corollary 3.2, it has been proved that this limit law is that of a
strong Markov process with semigroup $(p^\uparrow_t)$. Moreover,
from \cite{mi}, the process $(X,\Bbb{P}_x^{\mbox{{\bf
e}}/\varepsilon})$ leaves its minimum continuously, (that is
$\Bbb{P}_x^{ \mbox{{\bf e}}/\varepsilon}(X_m>X_{m-})=0$) if and
only if 0 is regular upwards. Then we  conclude using Proposition
\ref{prop1}.\QED

\noindent When $(X,\p)$ has no negative jumps and $0$ is not
regular upwards, the initial law of $(X,\Bbb{P}^{\uparrow })$ has
been computed in \cite{ch}. It is given by:
\begin{equation}
\Bbb{P}^{\uparrow }(X_{0}\in dx)=\frac{x\,\pi
(dx)}{\int_{0}^{\infty }u\,\pi (du)},\;\;\;x\ge 0\,,  \label{saut}
\end{equation}
where $\pi$ is the L\'evy measure of $(X,\p)$. It seems more
difficult to obtain an explicit formula which only involves $\pi$
in the general case.

We can now state our convergence result. Recall that it has been
proved in \cite{c1} in the special cases where $(X,\Bbb{P})$ is
either a L\'evy process which creeps downwards (see the next
section for the definition of the creeping of a stochastic
process) or a stable process. Note also that from Bertoin \cite
{b1}, Th. VII.14 this convergence also holds when $(X,\Bbb{P})$
has no positive jumps. Then Tanaka \cite{ta}, Theorems 4 and 5
proved the finite dimensional convergence in the very general
case.

\begin{theorem}
\label{main} Assume that $0$ is regular upwards. Then the family
$(\Bbb{P}_{x}^{\uparrow },x>0)$ converges on the Skorokhod space
to $ \Bbb{P}^{\uparrow }$. Moreover the semigroup
$(p_{t}^{\uparrow },t\geq 0)$ has the Feller property on
$[0,\infty )$.

If $0$ is not regular upwards, then for any $\varepsilon >0$ , the
process $(X\circ \theta _{\varepsilon },\Bbb{P}_{x}^{\uparrow })$
converges weakly towards $(X\circ \theta _{\varepsilon
},\Bbb{P}^{\uparrow }) $, as $x$ tends to $0$. In that case, the
Feller property of $ (p_{t}^{\uparrow },t\geq 0)$ holds only on
$(0,\infty )$.
\end{theorem}

\noindent {\it Proof}: Let $(\Omega,\mathcal{F},P)$ be a
probability space on which we can define a family of processes $
(Y^{(x)})_{x>0}$ such that each process $Y^{(x)}$ has law $\Bbb{P}
_{x}^\uparrow$. Let also $Z$ be a process with law
$\Bbb{P}^\uparrow$ and which is independent of the family
$(Y^{(x)})$. Let $m_x$ be the unique hitting time of the minimum
of $Y^{(x)}$ and define, for all $x>0$, the process $Z^{(x)}$ by:
\[
Z_t^{(x)}=\left\{
\begin{array}{ll}
Y_t^{(x)} & t<m_x \\
Z_{t-m_x}+Y_{m_x}^{(x)}\;\;\; & t\geq m_x\,.
\end{array}
\right.
\]
By the preceding theorem, under $P$, $Z^{(x)}$ has law
$\Bbb{P}_{x}^\uparrow$.

Now first assume that 0 is regular upwards, so that
$\lim_{t\downarrow 0}Z_{t}=0$, almost surely. We are going to show
that the family of processes $Z^{(x)}$ converges in probability
towards the process $Z $ as $x\downarrow 0$ for the norm of the
$J_{1}$-Skorohod topology on the space $\mathcal{D}([0,1])$. Let
$(x_{n})$ be a decreasing sequence of real numbers which tends to
0. For $\omega \in \mathcal{D}([0,1])$, we easily see that the
path $Z^{(x_{n})}(\omega )$ tends to $Z(\omega )$ as $n$ goes to
$\infty$ in the Skohorod's topology, if both $m_{x_{n}}(\omega )$
and $ \overline{Z}^{(x_{n})}_{m_{x_n}}(\omega )$ tend to 0. Hence,
it suffices to prove that both $m_{x}$ and
$\overline{Z}^{(x)}_{m_x}$ converges in probability to 0 as $
x\rightarrow 0$, i.e. for any fixed $\varepsilon >0,\eta >0,$
\begin{equation}
\lim_{x\downarrow 0}\Bbb{P}_{x}^{\uparrow }(m>\varepsilon
)=0\;\;\;\mbox{and}\;\;\;\mbox{ } \lim_{x\downarrow
0}\Bbb{P}_{x}^{\uparrow }(\overline{X}_{m}>\eta )=0. \label{limp}
\end{equation}
First, applying the Markov property at time $\varepsilon $ gives
\begin{eqnarray*}
\Bbb{P}_{x}^{\uparrow }(m >\varepsilon )&=&\int_{0<y\le
x}\int_{z>y}\Bbb{P} _{x}^{\uparrow }(X_{\varepsilon }\in
dz,\underline{X}_{\varepsilon }\in
dy,\varepsilon <\zeta )\Bbb{P}_{z}^{\uparrow }(U<y) \\
&=&\int_{0<y\le x}\int_{z>y}\Bbb{Q}_{x}(X_{\varepsilon }\in
dz,\underline{X} _{\varepsilon }\in dy,\varepsilon <\zeta
)\frac{h(z)}{h(x)}\Bbb{P}
_{z}^{\uparrow }(U<y) \\
&=&\int_{0<y\le x}\int_{z>y}\Bbb{P}_{x}(X_{\varepsilon }\in
dz,\underline{X} _{\varepsilon }\in dy)\frac{h(z)-h(z-y)}{h(x)},
\end{eqnarray*}
where we have used the result of Theorem \ref{theo1} and the fact
that $\Bbb{ Q}_{x}$ and $\Bbb{P}_{x}$ agree on
$\mathcal{F}_{\varepsilon }\cap (\underline{X}_{\varepsilon }>0)$.
Put $\tilde{h}:=h-h(0)$, and recall from Section 2 that
$\tilde{h}$ is increasing and subadditive, hence we have
$h(z)-h(z-y)\leq \tilde{h}(y),$ and so
\begin{eqnarray*}
\Bbb{P}_{x}^{\uparrow }(m &>&\varepsilon )\leq
\frac{1}{h(x)}\int_{0<y\le x}\int_{z>y}\Bbb{P}_{x}(X_{\varepsilon
}\in dz,\underline{X}_{\varepsilon
}\in dy)\tilde{h}(y) \\
&=&\frac{1}{h(x)}\int_{0<y\le
x}\Bbb{P}_{x}(\underline{X}_{\varepsilon }\in dy)\tilde{h}(y)\leq
\frac{\tilde{h}(x)}{h(x)}\Bbb{P}_{x}(\underline{X} _{\varepsilon
}>0)\,.
\end{eqnarray*}
When 0 is not regular downwards $h(0)=1$ (see Section 2), hence
$\frac{ \tilde{h}(x)}{h(x)}\rightarrow 0$ as $x\rightarrow 0$, and
we obtain the result in that case. When 0 is regular,
$\frac{\tilde{h}(x)}{h(x)}=1$, but in that case, we clearly have
$\Bbb{P}_{x}(\underline{X}_{\varepsilon }>0)\rightarrow 0$ as
$x\rightarrow 0$, so the result is also true.

For the second claim, we apply the strong Markov property at time
$\tau :=\tau _{(\eta ,\infty )}$, with $x<\eta$, to get
\begin{eqnarray*}
\Bbb{P}_{x}^{\uparrow }(\overline{X}_{m} &>&\eta )=\int_{z\geq
\eta }\int_{0<y\leq x}\Bbb{P}_{x}^{\uparrow }(X_{\tau }\in
dz,\underline{X}_{\tau
}\in dy,\tau<\zeta)\Bbb{P}_{z}^{\uparrow }(U<y) \\
&=&\int_{z\geq \eta }\int_{0<y\leq x}\Bbb{P}_{x}^{\uparrow
}(X_{\tau }\in dz, \underline{X}_{\tau }\in
dy,\tau<\zeta)\frac{h(z)-h(z-y)}{h(z)}.
\end{eqnarray*}
We now apply the simple bound
\[
\frac{h(z)-h(z-y)}{h(z)}\leq \frac{\tilde{h}(y)}{h(z)}\leq
\frac{\tilde{h}(x) }{h(\eta )}\mbox{ for }0<y\leq x\mbox{ and
}z\geq \eta
\]
to deduce that
\[
\Bbb{P}_{x}^{\uparrow }(\overline{X}_{m}>\eta
)\le\frac{\tilde{h}(x)}{h(\eta
)}\rightarrow0\;\mbox{as}\;x\downarrow0\,.
\]
Finally when 0 is not regular upwards, (\ref{limp}) still holds
but we can check that, at time $t=0$, the family of processes
$Z^{(x)}$ does not converge in probability towards 0. However
following the above arguments we can still prove that for any
$\varepsilon>0$, $(Z^{(x)}\circ\theta_\varepsilon)$ converges in
probability towards $Z\circ\theta_\varepsilon$ as $x\downarrow0$.
\QED

\noindent The following absolute continuity relation between the
measure $\underline{n}$ of the process of the excursions away from
0 of $ X-\underline{X}$ and $\Bbb{P}^{\uparrow }$ has been shown
in \cite{c1}: for $ t>0$ and $A\in \mathcal{F}_{t}$
\begin{equation}\label{eq2}
\underline{n}(A,t<\zeta )=k\Bbb{E}^{\uparrow }(h(X_{t})^{-1}A),
\end{equation}
where $k>0$ is a constant which depends only on the normalization
of the local time $\underline{L}.$ Relation (\ref{eq2}) was
established under some additional  hypotheses in \cite{c1} but we
can easily check that it still holds in all the cases which
concern us. Then a consequence of Theorem \ref {main} is:

\begin{corollary}
\label{coro1} For any $t>0$ and for any $\mathcal{F}_{t}$
-measurable, continuous and bounded functional $F$,
\[
\underline{n}(F,t<\zeta )=k\lim_{x\rightarrow
0}\Bbb{E}_{x}^{\uparrow }(h(X_{t})^{-1}F).
\]
\end{corollary}

\noindent Another application of Theorem \ref{main} is to the
asymptotic behavior of the semigroup $q_t(x,y)$, $t>0$, $y>0$,
when $x$ goes towards 0. Let us denote by $j_t(x)$, $t\ge0$,
$x\ge0$ the density of the entrance law of the excursion measure
$\underline{n}$, that is the Borel function which is defined for
any $t\ge0$ as follows:
\[
\underline{n}(f(X_t),t<\zeta )=\int_0^\infty f(x)j_t(x)\,dx\,,
\]
where $f$ is any positive or bounded Borel function $f$.

\begin{corollary}
\label{coro2} The asymptotic behavior of $q_{t}(x,y)$ is given by:
\[
\int_{0}^{\infty }f(y)q_{t}(x,y)\,dy\sim _{x\rightarrow
0}h(x)\int_{0}^{\infty }f(y)j_{t}(y)\,dy\,,
\]
for every continuous and bounded function $f$.
\end{corollary}

\noindent Note that when 0 is not regular downwards, the measure
$\underline{n}$ is nothing but $\Bbb{Q}_{0}$ and $h(0)=1$, so in
that case Corollaries \ref{coro1} and \ref{coro2} are
straightforward. In the other case, they are direct consequences
of Theorem \ref{main}, (\ref {eq1}) and (\ref{eq2}), so their
proofs are omitted.\\

\noindent We end this paper by a study of the asymptotic behaviour
of the function $h$ at 0 in terms of the lower tail of the height
of the generic reflected excursion. We define the height of a path
$\omega$ with finite lifetime $\zeta(\omega)$ as follows:
\[H(\omega):=\sup_{0\le t\le\zeta}\omega_t\,.\] The equality $n(H>x)h(x)=1$
is proved in \cite{b1}, Proposition VII.15 in the spectrally
negative case, i.e. when $(X,\p)$ has no positive jumps. However,
this relation does not hold in general; a counter example is
provided by the spectrally positive case, as we show hereafter: we
can prove, using Corollary \ref{coro1} that in any case one has
\begin{eqnarray*}
n(H>x)=n(\tau_{[x,\infty)}<\infty)&=&\lim_{y\downarrow0}\frac{1}{h(y)}
\Bbb{Q}_y(\tau_{[x,\infty)}<\tau_{(-\infty,0]})\nonumber\\
&=&\lim_{y\downarrow0}\frac{1}{h(y)}\p_y(\overline{X}(\tau_{(-\infty,0]})>x)\label{125}\,.
\end{eqnarray*}
Moreover, Theorem VII.8 in \cite{b1} implies that when the process
has no negative jumps,
\[
\p_y(\overline{X}(\tau_{(-\infty,0]})>x)=\frac{\hat{h}(x+y)}{\hat{h}(x)}\,,
\]
where $\hat{h}$ is the harmonic function defined as in Section 2
with respect to the dual process $\hat{X}\eqdef -X$. Then in this
 case, one has $h(y)=y$, so that from above,
$n(H>x)=\hat{h}'(x)/\hat{h}(x)$. We conclude that the expression
$n(H>x)h(x)$ can be equal to a constant only when $\hat{h}$ is of
the form $\hat{h}(x)=cx^\gamma$, for some positive constants $c$
and $\gamma$, but this is possible only if $(X,\p)$ is stable.
Note that the equality $n(H>x)h(x)=1$ has recently been noticed
independently  by Rivero \cite{ri} for any stable process, in a
work on the more general setting of Markov self-similar processes.

In the spectrally positive case discussed above, we can check that
the expression $n(H>x)h(x)$ tends to a constant as $x$ goes to 0
if moreover the process is in the domain of attraction of a stable
process, i.e. if there exists $\alpha\in(0,2]$ such that
$((t^{-1/\alpha}X_{st},\,s\ge0),\p)$ converges weakly to a stable
process with index $\alpha$, as $t$ goes to 0. This raised the
question of finding some other conditions under which this
property holds. The creeping of $(X,\p)$ is one such condition.

In the rest of the paper, we suppose that $(X,\p)$ does not drift
towards $-\infty$. We say that the process $(X,\p)$ creeps
$($upwards$)$ across the level $x>0$ if
\begin{equation}  \label{creep}
\p(X(\tau_{[x,\infty)})=x)>0\,.
\end{equation}
It is well known that if $X$ creeps across a positive level $x>0$,
then it creeps across all positive levels. Moreover this can
happen if and only if
\begin{equation}\label{id457}
\lim_{x\downarrow0}\p(X(\tau_{[x,\infty )})=x)=1\,.
\end{equation}
We refer for instance to  \cite{b1}, chap. VI for a proof of  this
equivalence. We need the following lemma for the proof of the next
theorem.

\begin{lemma}\label{lem45}
Define for all $x>0$,  $\sigma_x=\sup\{t:X_t\le x\}$. If $(X,\p)$
creeps upwards, then $\lim_{x\downarrow0}\Bbb{P}^{\uparrow
}(X(\sigma _{x})=x)=1$.
\end{lemma}

\noindent {\it Proof}: Let us first mention the following identity
which is derived from Theorems 4.1 and 4.2 of Duquesne \cite{du}:
\begin{equation}\label{id456}
(X(\sigma _{x}),\Bbb{P}^{\uparrow })=(X(\tau_{[x,\infty
)})+[X(g_{x})-X(\tau_{[x,\infty )}-)],\Bbb{P})\,,
\end{equation}
where $g_{x}=\sup \{t\le \tau_{[x,\infty
)}:\overline{X}_{t}=X_{t}\vee X_{t-}\}$. This identity can also be
checked through Doney-Tanaka's construction of L\'{e}vy processes
conditioned to stay positive (see Doney \cite{do}). Then observe
that $X(\tau_{[x,\infty )})\ge x$ and $X(g_{x})-X(\tau_{[x,\infty
)}-)\ge0$, $\p$-a.s. Moreover since $(X,\p)$ creeps upwards, it is
not a compound Poisson process, hence it cannot reach a positive
level for the first time by a jump. In particular, on the event
$\{X(\tau_{[x,\infty )})=x\}$, the process is $\p$-a.s. continuous
at time $\tau_{[x,\infty )}$ and on this event,
$X(g_{x})=X(\tau_{[x,\infty )}-)=X(\tau_{[x,\infty )})$ so that
\[\{X(\tau_{[x,\infty )})=x\}=\{X(\tau_{[x,\infty )})+[X(g_{x})-X(\tau_{[x,\infty
)}-)]=x\}\,,\;\;\;\mbox{$\p$-a.s.}\] We conclude with
(\ref{id456}) and (\ref{id457}). \QED

\begin{theorem}\label{theo123}
If $(X,\p)$ creeps upwards, then $\underline{n}
(H>x)h(x)\rightarrow 1$, as $x\rightarrow 0$.
\end{theorem}

\noindent {\it Proof}: Recall that under the hypothesis of the
theorem, the lifetime of $(X,\Bbb{P}^{\uparrow })$ is almost
surely infinite. Then fix $ x>0$. Since $\{H>x\}=\{\tau_{(x,\infty
)}<\zeta \}$, from the identity (\ref {eq2}) applied at the
stopping time $\tau_{(x,\infty)}$, we have
\begin{equation}
h(x)\,\underline{n}(H>x)=h(x)\,\underline{n}(\tau_{(x,\infty)}<\zeta)\le
\underline{n}(h(X(\tau_{(x,\infty )})),\,\tau_{(x,\infty)}<\zeta
)=1\,. \label{in4}
\end{equation}
From the Markov property applied at time $\tau_{(0,x]}$, and since
$h$ is increasing, we have for any $ 0<x\le y$,
\begin{eqnarray}
\Bbb{P}_{y}^{\uparrow }(X(\sigma _{x})=x) &=&\Bbb{E}_{y}^{\uparrow
}\left( \Bbb{P}_{X(\tau_{(0,x]})}^{\uparrow }(X_{\sigma
_{x}}=x)\mbox{\rm
1\hspace{-0.033in}I}_{\{\tau_{(0,x]}<\infty \}}\right)\nonumber \\
&=&\frac{1}{h(y)}\Bbb{E}_{y}^{\Bbb{Q}}\left(
h(X(\tau_{(0,x]}))\Bbb{P} _{X(\tau_{(0,x]})}^{\uparrow }(X_{\sigma
_{x}}=x)\mbox{\rm 1\hspace{-0.033in}I}
_{\{\tau_{(0,x]}<\zeta \}}\right)   \nonumber \\
&\le &\frac{h(x)}{h(y)}\,.  \label{in1}
\end{eqnarray}
On the other hand, from the Markov property at time
$\tau_{(x,\infty )}$ and (\ref{eq2}), we have under
$\Bbb{P}^{\uparrow }$:
\begin{eqnarray}
\Bbb{P}^{\uparrow }(X(\sigma _{x})=x)
&=&\Bbb{E}^{\uparrow}\left(\Bbb{P}
_{X(\tau_{(x,\infty )})}^{\uparrow }(X_{\sigma _{x}}=x)\right)\nonumber \\
&=&\underline{n}\left(h(X(\tau_{(x,\infty)}))
\Bbb{P}_{X(\tau_{(x,\infty)})}^{\uparrow }(X_{\sigma
_{x}}=x)\mbox{\rm 1\hspace{-0.033in}I} _{\{\tau_{(x,\infty
)}<\zeta \}}\right) \,.  \label{in2}
\end{eqnarray}
But since $X(\tau_{(x,\infty)})\ge x$, inequality (\ref{in1})
gives
\[
\Bbb{P}_{X(\tau_{(x,\infty)})}^{\uparrow}(X_{\sigma _{x}}=x)\le
h(x)/h(X(\tau_{(x,\infty)})),\ a.s.
\]
Hence, from (\ref{in2}) we have
\begin{equation}
\Bbb{P}^{\uparrow }(X(\sigma _{x})=x)\le
h(x)\,\underline{n}(\tau_{(x,\infty )}<\zeta )\,.\label{in3}
\end{equation}
Finally, we deduce the result from (\ref{in4}), (\ref{in3}) and
Lemma \ref{lem45}. \QED

\end{document}